\documentclass{article}
\usepackage{latexsym}
\usepackage{amsthm}
\usepackage{amssymb}
\newcommand{\cl}{C \kern -0.1em \ell}     
\newcommand{\mysim}{\kern -0.25em \sim}

 \newtheorem{definition}{Definition}
 
\title{Automorphic Forms and Dirac Operators on Conformally Flat Manifolds}

\author{R.S.~Krau{\ss}har \thanks{Lehrgebiet f\"ur Mathematik und ihre Didaktik, Erziehungswissenschaftliche Fakult\"at, Universt\"at Erfurt, Nordh\"auser Str. 63, D-99089 Erfurt, Germany.  E-mail: {\tt soeren.krausshar@uni-erfurt.de}}}

\begin{document}
\maketitle
\begin{abstract}
In this paper we present a summarizing description of the connection between Dirac operators on conformally flat manifolds and automorphic forms based on a series of joint work with John Ryan over the last fifteen years. We also outline applications to boundary value problems.  
\end{abstract}

\par\medskip\par
To Professor John Ryan for his 60th birthday for the long collaboration and friendship
\section{Introduction}    

A natural generalization to $\mathbb{R}^{n}$ of the classical
Cauchy-Riemann operator has proved to be the Euclidean Dirac
operator $D$. Here $\mathbb{R}^{n}$ is considered as embedded in the real $2^{n}$-dimensional  Clifford algebra $Cl_{n}$ satisfying the relation $x^{2}=-\|x\|^{2}$ for each $x\in\mathbb{R}^{n}$. The elements $e_{1},\ldots,e_{n}$ of the standard orthonormal basis of $\mathbb{R}^{n}$ satisfy the relation $e_{i}e_{j}+e_{j}e_{i}=-2\delta_{ij}$. The Dirac operator is defined to be $\sum_{j=1}^{n}e_{j}\frac{\partial}{\partial x_{j}}$. Clifford algebra valued functions $f$ and $g$ that satisfy $Df =0$ respectively $gD = 0$
are often called left (right) monogenic functions. 

Its associated function theory together with its
applications is known as Clifford analysis and can be regarded as
a higher dimensional generalization of complex function theory in
the sense of the Riemann approach. Indeed, associated to this
operator there is a higher dimensional direct analogue of Cauchy's
integral formula and other nice analogues, cf. \cite{dss}. As in complex analysis, also the Euclidean Dirac
factorizes the higher dimensional Euclidean Laplacian viz $D^2 = -
\Delta$. Indeed, the Euclidean Dirac operator has been used successfully in
understanding boundary value problems and aspects of classical
harmonic analysis in $\mathbb{R}^{n}$. See for instance
\cite{GS2,SW}.

\ On the other hand Dirac operators have proved to be extremely useful tools in
understanding geometry over spin and pin manifolds. 
Basic aspects of Clifford analysis over spin manifolds have been developed in \cite{ca,cn}.
Further in \cite{KraRyan1,KraRyan2,KraRyan3,KraJMAA} and elsewhere it is
illustrated that the context of conformally flat manifolds provide a useful setting for developing Clifford analysis.

\ Conformally flat manifolds are those manifolds which possess an
atlas whose transition functions are M\"obius transformations.
Under this viewpoint conformally flat manifolds can be regarded as
higher dimensional generalizations of Riemann surfaces.

Following the classical work of  N. H. Kuiper \cite{Kuiper}, one can construct examples of conformally flat manifolds by
factoring out a subdomain $U \subseteq \mathbb{R}^n$ by a torsion-free Kleinian
group $\Gamma$ acting totally discontinuously on $U$.

\par\medskip\par 
Examples of conformally flat manifolds include spheres, hyperbolas, real projective space, cylinders, tori, the M\"obius strip, the Kleinian bottle and
the Hopf manifolds $S^1 \times S^{n-1}$. The oriented manifolds among them are also
spin manifolds. In \cite{KraRyan1,KraRyan2,KraJMAA} explicit Clifford
analysis techniques, including Cauchy and Green type integral
formulas, have been developed for these manifolds.

\par\medskip\par 

Finally, in one of our follow-up papers \cite{BCKR} we also looked at a class of hyperbolic manifolds namely those that arise from factoring out upper
half-space in $\mathbb{R}^n$ by a torsion-free congruence subgroup, $H$, of the
generalized modular group $\Gamma_p$. $\Gamma_p$ is the arithmetic
group that is generated by $p$ translation matrices ($p < n$) and
the inversion matrix. In two real variables these are $k$-handled spheres. Notice that the group $\Gamma_p$ is not torsion-free, as it contains the negative identity matrix. Consequently, the topological quotient of upper half-space with $\Gamma_p$ has only the structure of an orbifold. To overcome this problem we deal with congruence subgroups of level $N \ge 2$, which are going to be introduced later on.
\par\medskip\par

In this paper we present an overview about some of our most important  joint results. In the final section of this paper we also outline some applications addressing boundary value problems modelling stationary flow problems on these classes of manifolds where we adapt the techniques from \cite{GS2} to this more general geometric context.

\section{Clifford algebras and spin geometry}
\subsection{Clifford algebras and orthogonal transformations}
As mentioned in the introduction,  we embed the $\mathbb{R}^n$ into the real Clifford algebra $Cl_n$ generated by the relation $x^2=-\|x\|^2$. For details, see \cite{cn,dss,GS2}. This relation defines the 
multiplication rules $e_i^2=-1$, $i=1,...,n$ and $e_i e_j = - e_j e_i$ $\forall i \neq j$. A vector space basis for $Cl_n$ is given by $1,e_1,\ldots, e_n, e_1 e_2,\ldots, e_{n-1}e_n,\ldots,$ $e_1\cdots e_n$.
Each $x \in \mathbb{R}^n \backslash\{0\}$ has an inverse of the form $x^{-1} =- \frac{x}{\|x\|^2}$.   
We also consider the reversion anti-automorphism defined by $\tilde{ab}=\tilde{b}\tilde{a}$, where $\tilde{e_j}=e_j$ $\forall j=1,...,n$ and the conjugation defined by
$\overline{ab} = \overline{b}\;\overline{a}$, where $\overline{e_j}=-e_j$  $\forall j=1,...,n$.
\par\medskip\par
Notice that $e_1 x e_1 = - x_1 e_1 + x_2 e_2 + \cdots + x_n e_n$. The multiplication of $e_1$ from the left and from the right realizes in a simple form a reflection in the $e_1$-direction. More generally, one can say: If $O \in O(n)$, then there are reflections  $R_1,\ldots,R_m$ such that $O=R_1 \cdots R_m$. In turn for each $R_j$ there exists a $y_j \in S^{n-1}$ such that $R_j x = y_j x y_j$ for all $x \in \mathbb{R}^n$. Summarizing, one can represent a general transformation of $O(n)$ in the way  
$O x = y_1 \cdots y_m x y_m \cdots y_1$, so $Ox=ax\tilde{a}$ with $a=y_1 \cdots y_m$.

This motivates the definition of the pin group as 
$$Pin(n+1)=\{a \in Cl_n \mid a = y_1 \cdots y_m,\; y_i \in S^n\}.$$

Each transformation of $O(n)$ can be written as  $Ox =ax\tilde{a}$ with an $a \in Pin(n)$.  In view of $ax\tilde{a}=(-a)x(-\tilde{a})$, $Pin(n)$ is a double cover of $O(n)$. A subgroup of index $2$ is the spin group defined by 
$$Spin(n) := \{a=y_1\cdots y_m \in Pin(n) \mid m \equiv 0(2)\}.$$

Again, $Spin(n)$ is the double cover of $SO(n)$. 

\subsection{Spin geometry}
Here we summarize some basic results from \cite{ca,ERy2007,Kuiper}. 

Let $M$ be a connected orientable Riemannian manifold with Riemann metric $g_{ij}$. 
 
Consider for $x \in M$ all orthonormal-bases of the tangential space $TM_x$, which again are mapped to orthonormal-bases of $TM_x$ by the action of the $SO(n)$. This gives locally rise to a fiber bundle. 
 
Gluing together all these fiber bundles gives rise to a principal bundle $P$ over $M$ with a copy of $SO(n)$. 
\par\medskip\par
This naturally motivates the question whether it is possible to lift each fiber to $Spin(n)$ in a continuous way to obtain a new principal bundle $S$ that covers $P$. 

However, the ambiguity caused by the sign may give a problem.  If $s:U \to  U \times SO(n)$ is a section then there are two options of lifting $s$ to a spinor bundle $s^*: U \to U \times Spin(n)$, namely $s^*$ and $-s^*$. So, it may happen that:
\begin{itemize}
\item
It is not always possible to choose the sign in order to construct in a unique way a bundle $S$ over $M$, such that each fiber is a copy of $Spin(n)$.
\item
There also might be several possibilities.
\item
The different spin structures are described by the cohomology group\\ $H^1(M,\mathbb{Z}_2)$. 
\subsection{The Atiyah-Singer-Dirac operator}
Let $M$ be a Riemannian spin manifold. Let $\Gamma$ be the Levi-Civita connection. Then $\Gamma g_{ij} = 0$. Stokes's theorem tells us that 
\begin{eqnarray*}
& & \int\limits_{\partial V} \langle s_1(x), n(x) s_2(x)\rangle_S d\sigma(x)\\ &=& \int\limits_{V} (\langle s_1(x) D,  s_2(x)\rangle_S + \langle s_1(x), D s_2(x)\rangle_S) dV
\end{eqnarray*}
The arising differential operator here is the Atiyah-Singer-Dirac operator. In a  local orthonormal basis $e_1(x),...e_n(x)$ the latter has the form 
$$
D = \sum\limits_{j=1}^n e_j(x) \Gamma^{*}_{e_j(x)}.
$$ 
\end{itemize}    

\subsection{The Dirac operator in $\mathbb{R}^n$ and $\mathbb{R} \oplus\mathbb{R}^{n}$}
Following \cite{dss} and others, the Dirac operator in $\mathbb{R}^n$ has the simple form 
$D = \sum\limits_{j=1}^n \frac{\partial }{\partial x_j} e_j$. 
In the so-called space of paravectors $\mathbb{R} \oplus \mathbb{R}^n$ it particularly has the form  
$D = \frac{\partial }{\partial x_0} + \sum\limits_{j=1}^n \frac{\partial }{\partial x_j} e_j$. The latter naturally generalizes the  well-known Cauchy-Riemann operator $\frac{\partial }{\partial x_0} + \frac{\partial }{\partial x_1}i$ in a straight forward way to higher dimensions by adding the additional basis elements. 

As mentioned in the introduction, functions $f: U \to Cl_n$ (where $U \subseteq \mathbb{R}^n$ resp. $U \subseteq \mathbb{R} \oplus \mathbb{R}^n$) that are in the kernel of $D$ are often called monogenic.

\section{Analysis on manifolds}

\subsection{Conformally flat manifolds in $\mathbb{R}^{n+1}$}

Following for example \cite{Kuiper}, conformally flat manifolds are Riemannian manifolds with vanishing Weyl tensor. In turn these are exactly those 
 Riemannian manifolds which have an atlas whose transition functions are conformal  maps. In 
$\mathbb{R}^2$ conformal maps are exactly the (anti-)holomorphic functions. So, one deals with classical Riemann surfaces. Up from space dimension $n \ge 3$ in  
$\mathbb{R}^{\ge 3}$ however, the set of conformal maps coincides with the set of M\"obius transformations, cf. \cite{cn} The latter set of functions then represent reflections at spheres and hyperplanes. This seems to be quite restrictive at the first glance. However, this is not the case as the abundance of classical important examples will show as outlined in the following.  
\par\medskip\par
So, let us now turn to an explicit construction principle of conformally flat manifolds in $\mathbb{R}^n$ with $n \ge 3$. 
To proceed in this direction take a torsion free discrete subgroup ${\mathcal{G}}$ of the orthogonal group $O(n+1)$ that acts totally discontinuously on a simply connected domain ${\mathcal{D}}$. Next define a group action ${\mathcal{G}} \times {\mathcal{D}} \to {\mathcal{D}}$. Then the topological  quotient space ${\mathcal{D}}/{\mathcal{G}}$ is a conformally flat manifold.
As also shown in the original work of N.H.~Kuiper in 1949 \cite{Kuiper}, the universal cover of a conformally flat manifold possesses a local diffeomorphism to $S^{n}$. 
Conformally flat manifolds of the form ${\mathcal{D}}/{\mathcal{G}}$ are exactly those for which this local diffeomorphism is a covering map ${\mathcal{D}} \rightarrow {\mathcal{D}} \subset S^n$. 
\par\medskip\par 
Let us present a few elementary examples: 
\begin{itemize}
\item Take ${\mathcal{G}}={\mathcal{T}}_{p} := \mathbb{Z} + \mathbb{Z} e_1 +\cdots + \mathbb{Z} e_{p-1}$,   ${\mathcal{D}} = \mathbb{R}^{n+1}$ and consider the action
$$
(m_0,\ldots,m_n) \circ (x_0,\ldots,x_n) \mapsto (x_0+m_0,\ldots,x_{p-1}+m_{p-1},x_p,\ldots,x_n).
$$
Then the topological quotients ${\mathcal{D}}/{\mathcal{G}}$ represent oriented $p$ cylinders. These are spin manifolds with $2^{p}$ many different spinor bundles.  
In the particular case $p=n+1$ one deals with a compact oriented $n+1$-torus, cf. \cite{KraRyan1,KraRyan2}.  
\item
Take again as group ${\mathcal{T}}_{n+1}$ and the same domain, but consider a different action of the form 
$$
(m_0,\ldots,m_{n-1},m_n) \circ (x_0,\ldots,x_{n-1},x_n)$$ 
$$\mapsto (x_0+m_0,\ldots,x_{n-1}+m_{n-1},(-1)^{m_n} x_n + m_n).
$$
Now ${\mathcal{D}}/G$ is the non-orientable Kleinian bottle in $\mathbb{R}^{n+1}$. Due to the lack of orientability it is not spin. However, it is a pin manifold with $2^{n+1}$ many  pinor bundles. For some Clifford analysis on these manifolds we refer the reader to our recent works \cite{KraJMAA,KRV}.
\end{itemize}

Further interesting examples can be constructed by forming quotients with non-abelian subgroups of M\"obius transformations in $\mathbb{R}^n$, in particular with arithmetic subgroups of the so-called Ahlfors-Vahlen group discussed for instance in \cite{EGM90} and many other papers. To make the paper self-contained we recall its definition:

\begin{definition} (Ahlfors-Vahlen group)\\
A Clifford algebra valued matrix 
$M = \left(\begin{array}{cc} a & b \\ c & d \end{array} \right) \in Mat(2,Cl_n)$ 
belongs to the special Ahlfors Vahlen group $SAV(n)$, if:\\
- $a,b,c,d$ are products of paravectors from $\mathbb{R} \oplus \mathbb{R}^n$\\
- $a\tilde{c},c\tilde{d},d\tilde{b},b\tilde{c} \in \mathbb{R} \oplus \mathbb{R}^n$ and $a\tilde{d}-b\tilde{c}=1$. 
\end{definition}

The use of Clifford algebras allow us to describe M\"obius transformations in $\mathbb{R} \oplus \mathbb{R}^n$ in the simple way $M<x> = (ax+b)(cx+d)^{-1}$, similarly to complex analysis. For our needs we consider the special hypercomplex modular group \cite{KraHabil} defined by 
$$
{\Gamma_p} := \Bigg\langle \underbrace{\left(\begin{array}{cc}
1 & 1 \\ 0 & 1\end{array}\right), \left(\begin{array}{cc} 1 & e_1
\\ 0 & 1\end{array}\right),\ldots,\left(\begin{array}{cc} 1 & e_p
\\ 0 &
1\end{array}\right)}_{=:\mathcal{T}_p},\underbrace{\left(\begin{array}{cc}
0 & -1 \\ 1 & 0\end{array}\right)}_{=:J}\Bigg\rangle
$$
If $p < n$, then by applying the same arguments as in \cite{Freitag}, $\Gamma_p$ acts totally discontinuously on upper half-space discussed in \cite{EGM90} 
$$
H^{+}(\mathbb{R} \oplus\mathbb{R}^n) :=\{x_0 + x_1 e_1 +\cdots + x_n e_n \in \mathbb{R} \oplus\mathbb{R}^n \mid x_n > 0\}.
$$
In two dimensions it coincides with the classical group $SL(2,\mathbb{Z})$. To get a larger amount of examples we look at the following arithmetic congruence subgroups:
$$
\Gamma_p[N] := \Bigg\{ \left(\begin{array}{cc} a & b \\ c & d \end{array} \right) \in \Gamma_p\;|\;a-1,b,c,d-1 \in N \mathcal{O}_p\Bigg\},$$
where 
$\mathcal{O}_p  :=  \sum\limits_{A \subseteq P(\{1,\ldots,p\})} \mathbb{Z} e_A$ are the standard orders in $Cl_n$. \par\medskip\par
If we now take ${\mathcal{D}} = H^+(\mathbb{R}^{n+1})$, $G=\Gamma_p[N]$ with $N \ge 2$ and if we consider the action    
$$
(M,x) \mapsto M<x> := (ax+b) \cdot (cx+d)^{-1},
$$
where $\cdot$ is the Clifford-multiplication, then for $N\ge 2$ the topological quotient ${\mathcal{D}}/{\mathcal{G}}$ realizes a class of conformally flat orientable manifold with spin structures generalizing $k$-tori to higher dimensions, cf. \cite{BCKR}.

\subsection{Sections on conformally flat manifolds}
In the sequel let us make the general assumption that ${\mathcal{M}}:={\mathcal{D}}/{\mathcal{G}}$ is an orientable conformally flat manifold. 
Let furthermore $f:{\mathcal{D}} \to \mathbb{R}^{n+1}$ be a function with  $f(G<x>))=j(G,x) f(x)$ $\forall G \in {\mathcal{G}}$ where $j(G,x)$ is an automorphic factor satisfying a certain co-cycle relation. In the simplest case $j(G,x) \equiv 1$ one deals with a totally invariant function under the group action of $G$. Then the 
canonical projection $p:{\mathcal{D}} \to  {\mathcal{M}}$ induces a well-defined function $f':=p(f)$ on the quotient manifold ${\mathcal{M}}$. Now let $D:=\frac{\partial }{\partial x_0} + \sum\limits_{i=1}^n \frac{\partial }{\partial x_i} e_i$ be the Euclidean Dirac- or Cauchy-Riemann operator and let $\Delta$ be the usual Euclidean Laplacian. The  
canonical projection $p:{\mathcal{D}} \to  {\mathcal{M}}$ in turn induces a Dirac operator $D'=p(D)$ resp. Laplace operator $\Delta'=p(D)$ on  ${\mathcal{M}}$.

{\it Consequence}: A ${\mathcal{G}}$-invariant function on ${\mathcal{D}}$ from Ker $D$ (resp. from Ker $\Delta$) induces  functions on ${\mathcal{M}}$ in Ker $D'$ resp. in Ker $\Delta'$.  
More generally, one considers for $f'$ monogenic resp. harmonic spinor sections with values in certain spinor bundles.

\section{Automorphic forms}
\subsection{Basic background}
Following classical literature on automorphic forms, for instance \cite{Freitag}, 
let ${\mathcal{G}}$ be a discrete group that acts  totally discontinuously on a domain ${\mathcal{D}}$ by ${\mathcal{G}}
\times {\mathcal{D}} \to {\mathcal{D}},\;(g,d) \to d^*$, Roughly speaking, 
automorphic forms on ${\mathcal{G}}$ are functions on    
${\mathcal{D}}$ that are  quasi-invariant under the action of ${\mathcal{G}}$. The fundamental theory of holomorphic automorphic forms in one complex variable was founded around 1890, mainly by H. Poincar\'e, F. Klein and R.
Fricke. The associated quotient manifolds are holomorphic Riemann surfaces. 
\par\medskip\par
The simplest examples are holomorphic periodic functions. They serve as example of automorphic functions on discrete translation groups. Further very classical examples are the so-called  
Eisenstein series $$G_n(\tau) := \sum\limits_{(c,d) \in
\mathbb{Z}^2 \backslash \{(0,0)\}} (c\tau +d)^{-n} \;\;\; n \equiv
0(\;{\rm mod}\;2)\;\;n \ge 4$$ 
They are holomorphic functions on    
$H^+(\mathbb{C}) := \{z = x+iy \in \mathbb{C} \mid y > 0\}$. For all $z \in H^+(\mathbb{C})$ they satisfy :
$$
f(z) = (f|_n M)(z) \;\;\forall M = \left(\begin{array}{cc} a & b
\\ c & d \end{array}\right) \in SL(2,\mathbb{Z})
$$
where $(f|_n M)(z) := (cz+d)^{-n}
f\Big(\frac{az+b}{cz+d}\Big)$.

The Eisenstein series $G_n$ are the simplest non-vanishing holomorphic automorphic forms on    
$SL(2,\mathbb{Z})$. Further important examples are Poincar\'e series: For $n > 2,\;n \in 2 \mathbb{N}$ let  
$$
P_n(z,w) = \sum\limits_{M \in SL(2,\mathbb{Z})}
(cz+d)^{-n}(w+M<z>)^{-n}
$$
These functions have the same transformation behavior as the previously described Eisenstein series, namely: 
$$
P_n(w,z) = P_n(z,w) = (cz+d)^{-n}P_n(M<z>,w).
$$
In contrast to the Eisenstein series, these Poincar\'e series have the special property that they vanish at the singular points of the quotient manifold resp. orbifold.  
\subsection{Higher dimensional generalizations in $n$ real variables}
Already in the 1930s C.L. Siegel developed higher dimensional analogues of the classical automorphic forms in the framework of holomorphic functions in several complex variables. The context is Siegel upper half-space where one considers the action of discrete subgroups of the symplectic group.
\par\medskip\par
More closely related to our intention is the line of generalization initiated by H. Maa{\ss} and extended by  J. Elstrodt, F. Grunewald, J.
Mennicke \cite{EGM90} and A. Krieg \cite{Kri88} and many followers in the period of 1985-1990 and onwards.\\ These authors looked at higher dimensional generalizations of 
Maa{\ss} forms which are non-holomorphic automorphic forms on discrete subgroups of the Ahlfors-Vahlen group (including for example the Picard group and the Hurwitz quaternions) on upper half-space $H^+(\mathbb{R}^n)$. As important analytic properties they exhibit to be complex-valued eigensolutions to the scalar-valued Laplace-Beltrami operator  
$$
\Delta_{LB} = x_n^2 \Big(\sum\limits_{i=0}^n
\frac{\partial^2}{\partial x_i^2}\Big) -(n-1) x_n \frac{\partial
}{\partial x_n}
$$
In the case $n=1$ one has $\Delta_{LB} = x_1^2 \Delta$. Therefore, in the one-dimensional case holomorphic automorphic forms simply represent a special case of  Maa{\ss} forms.
\par\medskip\par
A crucial disadvantage of  Maa{\ss} forms behind the background of our particular intentions consists in the fact that they 
do not lie in the kernel of the Euclidean Dirac or Laplace operator. Additionally, they are only scalar-valued.

\subsection{Clifford-analytic automorphic forms}

{\it Some milestones in the literature}\\[0.3cm]
To create a theory of Clifford algebra valued automorphic forms that are in kernels of Dirac operators remained a challenge for a long time. A serious obstacle has been to overcome the problem that neither the multiplication nor the composition of monogenic functions are monogenic again. However, the set of monogenic functions is quasi-invariant under the action of M\"obius transformations up to an automorphic factor that fortunately obeys a certain co-cyle relation. The latter actually provides the fundament to build up a theory of automorphic forms in the Clifford analysis setting. 
\par\medskip\par

The first contribution in the Clifford analysis setting is the paper \cite{Ry1983} where J. Ryan constructed $n$-dimensional monogenic analogues of the Weierstra{\ss} $\wp$-function and the Weierstra{\ss} $\zeta$-function built as summations of the monogenic Cauchy kernel over an $n$-dimensional lattice. Here, the invariance group is an abelian translation group with $n$ linearly independent generators. 
\par\medskip\par
In the period of 1998-2004 the author developed the fundamentals for a more general theory of Clifford holomorphic automorphic forms on more general arithmetic subgroups of the Ahlfors-Vahlen group, cf. \cite{KraHabil}. The geometric context is again upper half-space. However, the functions are in general Clifford algebra valued and are null-solutions to iterated Dirac equations. In fact, the framework of null-solutions to the classical first order Dirac opertor is too small for a large theory of automorphic forms, because the Dirac operator admits only the construction of automorphic forms to one weight factor only. To consider automorphic forms with several automorphic weight factors a more approriate framework is the context of iterated Dirac equations of the form $D^m f=0$ to higher order integer powers $m$. In fact, higher order Dirac equations can also be related to $k$-hypermonogenic functions \cite{Leutwiler} which allows us to relate the theory of Maa{\ss} forms to Clifford holomorphic automorphic forms. This is successfully explained in \cite{CKR2007,CGK2013}. In this context it was possible to generalize Eisenstein- and Poincar\'e series construction which gave rise to the construction of spinor sections with values in certain spinor bundles on the related quotient manifolds. 
\par\medskip\par
In fact, as shown in \cite{CGK2013}, it is possible to decompose the Clifford module of Clifford holomorphic automorphic forms as an orthogonal direct sum of Clifford holomorphic Eisenstein- and Poincar\'e series. As shown in our recent paper \cite{GK2015}, both modules turn out to be finitely generated.    
\par\medskip\par    
The applications to solve boundary value problems on related spin manifolds started with our first joint paper \cite{KraRyan1} where we applied multiperiodic Eisenstein series on translation groups to construct Cauchy and Green kernels on confomally flat cylinders and tori associated to the trivial spinor bundle. 
\par\medskip\par
In \cite{KraRyan2} we extended our study to address the other spinor bundles on these manifolds. Furthermore, we also looked at dilation groups instead of translation groups, too, and managed to give closed representation formulas for the Cauchy kernel of the Hopf manifold $S^1 \times S^{n-1}$. Already in this paper we outlined the construction of spinor sections on some hyperbolic manifolds of genus $\ge 2$. However, to obtain the Cauchy kernel for these manifolds we needed to construct Poincar\'e series which was a hard puzzle piece to find. Eisenstein series were easy to construct and they lead to non-trivial spinor sections on these manifolds. However, the Cauchy kernel must have the property to vanish at the cusps of the group - and that construction was hard to do. A breakthrough in that direction could be established in our joint paper \cite{BCKR} where we were able to fill in that gap. 
\par\medskip\par
Finally, in \cite{CKR2012} and \cite{CGK2013} we were able to extend these constructions to the framework of $k$-hypermonogenic functions \cite{slnew} and holomorphic Clifford functions addressing null-solutions to $D \Delta^m f = 0$. This is the function class considered by G. Laville and I. Ramadanoff introduced in \cite{LR}.        
\par\medskip\par
 Summarizing, the theory of Clifford holomorphic automorphic form provides us with a toolbox to solve boundary value problems related to the Euclidean Dirac or Laplace operator on conformally flat spin manifolds or more generally on manifolds generalizing classical modular curves.
\par\medskip\par
It also turned out to be possible to make analogous constructions for some non-orientable conformally flat manifolds with pin structures.
Belonging to that context, in \cite{KraRyan2} we addressed real projective spaces and in \cite{KraJMAA,KRV} we carried over our constructions to the framework of higher dimensional M\"obius strips and the Klein bottle. 
  
\subsection{Some concrete examples}
\begin{itemize}
\item
The simplest non-trivial examples of Clifford-holomorphic automorphic forms on the translation groups ${\mathcal{T}}_p$ are given by the series 
$$
\epsilon^{(p)}_{\bf m}(x):=\sum_{\omega \in \Omega_p} 
q_{\bf m}(x+\omega),\quad\quad \Omega_p:= \mathbb{Z} \omega_0 +
\mathbb{Z}\omega_1 + \ldots \mathbb{Z} \omega_p
$$
which are normally convergent, if $|{\bf m}| \ge p+2$. Here, $q_{\bf m}(x) := \frac{\partial^{|{\bf m}|}}{\partial x^{\bf m}} q_{\bf 0}(x)$ where $q_{\bf 0}(x) := \frac{\overline{x}}{\|x\|^{m+1}}$ and where ${\bf m} :=(m_1,\cdots,m_n)$ is a multi-index and $|{\bf m}| := m_1+m_2+\cdots+m_n$ and $x^{\bf m} := x_1^{m_1} \cdots x_n^{m_n}$ is used as in usual mult-index notation. 

These series $\epsilon^{(p)}_{\bf m}(x)$ can be interpreted as 
Clifford holomorphic generalizations of the trigonometric  functions (with singularities) and of the  
 Weierstra{\ss} $\wp$-function. The projection $p(\epsilon^{(p)}_{\bf m}(x))$ then induces a well-defined spinor section on the cylinder resp. torus $\mathbb{R}^{n+1}/\Omega_p$ with values in the trivial spinor bundle. Other spinor bundles can be constructed by introducing proper minus signs and the following decomposition of the period lattice in the way
$\Omega_p := \Omega_l \oplus \Omega_{p-l}$ where $0\le l \le p$.  The proper analogues of $\epsilon^{(p)}_{\bf m}(x)$ then are defined by
$$
\epsilon^{(p,l)}_{\bf m}(x):=\sum_{\omega \in \Omega_l \oplus \Omega_{p-l}}(-1)^{m_0+\cdots+m_l} 
q_{\bf m}(x+\omega),\quad \Omega_p:= \mathbb{Z} \omega_0 +
\mathbb{Z}\omega_1 + \ldots \mathbb{Z} \omega_p.
$$ 
In total one can construct $2^{p+1}$ different spinor bundles. The Cauchy kernel is obtained by the series that one obtains in the case ${\bf m}={\bf 0}$. In that case it may happen that the associated series above does not converge. To obtain convergence in those cases one needs to apply special convergence preserving terms. For the technical details we refer to our papers \cite{KraRyan1,KraRyan2,KraRyan3}.   

\item Let us now turn to examples of hyperbolic manifolds that are generated by quotient forming with non-abelian groups. The simplest non-trivial example in that context are the  monogenic Eisenstein series for $\Gamma_{n-1}[N]$ (and also for $\Gamma_{p}[N]$) with $p <n-1$ which have been introduced in \cite{CKR2007}.
$$
E(z)=\lim\limits_{\sigma\rightarrow 0^+} \sum\limits_{M:
\mathcal{T}_{n-1}[N]\backslash \Gamma_{n-1}[N]}
\Bigg(\frac{x_n}{\|cx+d\|^2}  \Bigg)^{\sigma} q_{\bf 0}(cx+d).
$$
Notice that we here applied the Hecke trick (cf. \cite{Freitag}) to get well-definedness. In fact, these Eisenstein series provide us with the simplest examples of non-vanishing spinor sections defined on the hyperbolic quotient manifolds $H^+(\mathbb{R}^{n+1})/{\Gamma_p[N]}$. However, these series do not vanish at the cusps of the group. They serve as examples but they don't reproduce the Cauchy integral. This property can be achieved by the following  monogenic Poincar\'e series, introduced in our papers \cite{CKR2007,BCKR}. The latter have the form 
$$
P_p(x,w) = \lim_{\sigma\rightarrow 0^+}\sum_{M \in
\Gamma_p[N]}\Bigg(\frac{x_n}{\|cx+d\|^2} \Bigg)^{\sigma}q_{\bf
0}(cx+d)q_{\bf 0}(w+M\langle x\rangle).
$$
The series $P_p(x,w)$ are indeed monogenic cusp forms, in particular  
$$\lim\limits_{x_n\rightarrow \infty}P_p(x_n e_n)=0.$$ Its projection down to the manifold induce the Cauchy kernel. In the following section we want to outline how the explicit knowledge of the Cauchy kernel allows us to solve boundary value problems on these manifolds. 
\end{itemize}
\section{Applications to BVP on manifolds}
Some practical motivations for the following studies are to understand weather forecast and flow problems on spheres, cylindrical ducts or other geometric contexts fitting into the line of investigation of \cite{Sprweather,CKK2015}. 
\subsection{The Cauchy kernel on spin manifolds}
Monogenic generalization of the Weierstra{\ss} $\wp$-function $\varepsilon^{(p)}_n$ give rise to monogenic sections on  $p$-cylinders   $\mathbb{R}^{n+1}/{\mathcal{T}}_p$. So, 
monogenic automorphic forms on  $\Gamma_p({\mathcal{I}})[N]$ define spinor sections on the  $k$-tori  $H^{+}(\mathbb{R}^{n+1})/\Gamma_p({\mathcal{I}})[N]$. The Poincar\'e series give us the Cauchy kernel on these manifolds. Summarizing, 
for $p < n-1$ the Cauchy kernel has the concrete and explicit form    
$$
C(x,y) = \sum\limits_{M \in \Gamma_p({\mathcal{I}})[N]} \frac{\overline{cx+d}}{|cx+d|^{n}} \frac{(\overline{y-M<x>})}{\|y-M<x>\|^{n+1}}.
$$
Each monogenic section $f'$ on ${\mathcal{M}}$ then satisfies  
$$
f'(y') = \int\limits_{\partial S'} C'(x',y')d\sigma'(x') f'(x')
$$
which is the reproduction of the Cauchy integral, cf. \cite{BCKR}.  
\subsection{The stationary Stokes flow problem on some conformally flat spin manifolds}

Let ${\mathcal{M}}:={\mathcal{D}}/G$ be a conformally flat spin manifold for which we know the Cauchy kernel $C(x,y)$ to the Dirac operator (concerning a fixed chosen spinor bundle $F$). 

Next let $E \subset {\mathcal{M}}$ be a domain with suffiently smooth boundary $\Gamma:=\partial E$.

Now we want to consider the following Stokes problem on  ${\mathcal{M}}$:
\begin{eqnarray}\label{Eq.1}
-\Delta u + \frac{1}{\eta} p &=& F \;\;{\rm in}\;E \\
{\rm div}\;u &=& 0 \;\;{\rm in}\;E\nonumber \\
u &=& 0 \;\;{\rm on}\;\Gamma, \nonumber
\end{eqnarray}
where $u \in W^{1}_2(E,F)$ is the velocity of the flow and where $p \in W^{1}_2(E,\mathbb{R})$ is the hydrostatic pressure. The    explicit knowledge of the Cauchy kernel on ${\mathcal{M}}$ allows us to set up explicit analytic representation formulas for the solutions.\par\medskip\par
To meet this end we define in close analogy of \cite{GS2} the Teodorescu transform and the Cauchy transform on ${\mathcal{M}}$ by     
\begin{eqnarray*}
(T_E f)(x) &:=&-\int_{E} C(x,y) f(y) dV(y)\\
(F_{\Gamma} f)(x) &:=& \int_{\Gamma} C(x,y) d\sigma(y) f(y).
\end{eqnarray*}
where now $C(x,y)$ stands for the Cauchy kernel associated to the chosen spinor bundle $F$ on the manifold. 
The associated Bergman projection $P:L^2(E) \to L^2(E)\cap Ker(D)$ then can be represented in the usual form 
$$
P=F_{\Gamma}(tr_{\Gamma}T_E F_{\Gamma})^{-1}tr_{\Gamma}T_{E}.
$$
An application of the Clifford analysis methods provide us with explicit analytic formulas for the velocity and the pressure of the fluid running over this manifold.An application of $T_E$ to (\ref{Eq.1}) yields:
$$
(T_E D)(Du) + \frac{1}{\eta} T_E D p = T_E F.
$$
Next, an application of the Borel-Pompeiu formula (Cauchy-Green formula) leads to:
$$
Du - F_{\Gamma} Du + \frac{1}{\eta} p - \frac{1}{\eta} F_{\Gamma} p = T_E F
$$
If we apply the Pompeiu projection $Q:=I-P$, then one obtains  
$$
QDu - Q F_{\Gamma} Du + \frac{1}{\eta} Q p - \frac{1}{\eta} Q F_{\Gamma} p = Q T_E F \;\;(*)
$$
Since $F_{\Gamma} Du, F_{\Gamma} p \in $ Ker $D$, one further gets   
$$
Q F_{\Gamma} Du = 0 \;\;\;{\rm and}\;\;\; Q F_{\Gamma} p = 0.
$$
Thus, \begin{equation} 
QDu + \frac{1}{\eta} Q p = Q T_E F.
\end{equation}
Further, another application of $T_E$ and the property $QDu= Du$  leads to  
\begin{eqnarray*}
  T_E Du + \frac{1}{\eta} T_E Q p &=& T_E Q T_E F\\
\Leftrightarrow u - \underbrace{F_{\Gamma} u}_{=0} + \frac{1}{\eta} T_E Q p &=& T_E Q T_E F.
\end{eqnarray*}
In view of $u|_{\Gamma} = 0$ we obtain the following representation formula for the velocity:     
$$
{u} = T_E (I - P) T_E F - \frac{1}{\eta} T_E(I-P) p.
$$
Finally, the condition div $u=0$ allows us to determinate the pressure 
$$
(\Re(I-P)) p = \eta \Re((I-P) T_E (I-P) F)
$$
where $\Re(\cdot)$ stands for the scalar part of an element from the Clifford algebra.
\par\medskip\par 
{\bf Remark}. An extension of this approach to the parabolic setting in which instationary flow problems are considered are treated in our follow up paper \cite{CKK2015}. 

\section{Acknowledgments}

The work of the author is supported by the project \textit{Neue funktionentheoretische Methoden f\"ur instation\"are PDE}, funded by Programme for Cooperation in Science between Portugal and Germany, DAAD-PPP Deutschland-Portugal, Ref: 57340281.

\end{document}